\documentclass{article}
\usepackage{amsmath}
\usepackage{amsfonts}
\usepackage{amssymb}
\usepackage{graphicx}
\providecommand{\U}[1]{\protect \rule{.1in}{.1in}}
\newtheorem{theorem}{Theorem}
\newtheorem{corollary}[theorem]{Corollary}

\newtheorem{proposition}[theorem]{Proposition}

\newenvironment{proof}[1][Proof]{\noindent \textbf{#1.} }{\  \rule{0.5em}{0.5em}}
\begin{document}

\title{On the general form of bimonotone operators}
\author{Nicolas Hadjisavvas\thanks{Department of Product and Systems Design
Engineering, University of the Aegean, Syros, Greece.
Email: nhad@aegean.gr}}
\date{}
\maketitle

\begin{abstract}
	In a recent paper (2024) Camacho, C\'{a}novas, Mart\'{\i}nez-Legaz and Parra introduced bimonotone operators, i.e., operators $T$ such that both $T$ and $-T$ are monotone, and found some interesting applications to convex feasibility problems, especially in the case the operator is also paramonotone. In the present paper we drop paramonotonicity and examine the question of finding the most general form of a bimonotone operator in a Banach space. We show that any such operator can be reduced in some sense to a single-valued, skew symmetric linear operator. This facilitates the proof of some results involving these operators in applications.
\end{abstract}

\section{Introduction and preliminaries}

Let $X$ be a normed space and $X^{\ast}$ be its dual. Given a multivalued
operator $T:X\rightrightarrows X^{\ast}$ we will denote by
$\operatorname*{dom}T$, $\operatorname*{range}T$ and $\operatorname*{gph}T$
its domain, range and graph, respectively. We recall that $T$ is called
monotone if for every $\left(  x,x^{\ast}\right)  $, $\left(  y,y^{\ast
}\right)  \in \operatorname*{gph}T$, one has $\left \langle x^{\ast}-y^{\ast
},x-y\right \rangle \geq0$. It is called bimonotone \cite{Martinez} if both $T$
and $-T$ are monotone, i.e.,
\begin{equation}
	\forall \left(  x,x^{\ast}\right)  ,\left(  y,y^{\ast}\right)  \in
	\operatorname*{gph}T,\quad \left \langle x^{\ast}-y^{\ast},x-y\right \rangle =0.
	\label{bimon}
\end{equation}

The operator $T$ is called paramonotone \cite{CensIusZen}, if the following
implication holds:
\begin{equation}
	\left.
	\begin{array}
		[c]{c}
		\left(  x,x^{\ast}\right)  ,\left(  y,y^{\ast}\right)  \in \operatorname*{gph}
		T\\
		\left \langle x^{\ast}-y^{\ast},x-y\right \rangle =0
	\end{array}
	\right \}  \Rightarrow \left(  x,y^{\ast}\right)  ,\left(  y,x^{\ast}\right)
	\in \operatorname*{gph}T. \label{paramon}
\end{equation}

In a recent paper \cite{Martinez} the authors noted that, as a consequence of
(\ref{bimon}) and (\ref{paramon}), an operator which is both bimonotone and
paramonotone is constant on its domain. This observation leads to an
interesting application to convex feasibility problems, regarding translations
of mutually disjoint convex sets.

In this paper we answer the following question: What if we drop the assumption of
paramonotonicity? What is the general form of an arbitrary bimonotone
operator? In the special case of a single-valued, bimonotone operator $T$
defined on the whole space $\mathbb{R}^{n}$, the answer is known: there exists
a linear, skew-symmetric operator $A$ on $\mathbb{R}^{n}$ and a vector
$v\in \mathbb{R}^{n}$ such that $T(x)=Ax+v$ for all $x\in \mathbb{R}^{n}$ (see
Prop. 3.1 in \cite{BiaHadScha}). In this paper, we answer the above question
in its most general form, for an arbitrary multivalued operator.

We now give some notation and definitions. If $V$ is a linear (not necessarily
closed) subspace of $X$, let $J$ be the embedding $J:V\rightarrow X$ which is
simply the operator $J:V\ni x\rightarrow x\in X$. Its adjoint $J^{\ast
}:X^{\ast}\rightarrow V^{\ast}$ is the operator such that $\left \langle
J^{\ast}x^{\ast},x\right \rangle =\left \langle x^{\ast},Jx\right \rangle $ for
every $x\in V$, $x^{\ast}\in X^{\ast}$, i.e., $J^{\ast}x^{\ast}$ is the
restriction of $x^{\ast}$ onto $V$.

Given an operator $T:X\rightrightarrows X^{\ast}$ with $0\in
\operatorname*{dom}T$, let $V=\left[  \operatorname*{dom}T\right]  $ be the
linear span of $\operatorname*{dom}T$. The operator $\hat{T}
:V\rightrightarrows V^{\ast}$ such that $\hat{T}\left(  x\right)  =J^{\ast
}T\left(  Jx\right)  =\left \{  J^{\ast}x^{\ast}:x^{\ast}\in T\left(
Jx\right)  \right \}  $ will be called the reduction of $T$ on $V$. That is,
$\hat{T}\left(  x\right)  $ is the set of all restrictions on $V$ of elements
of $T\left(  x\right)  $.

A linear operator $A:V\rightarrow V^{\ast}$ is called skew symmetric if
$\left \langle Ax,y\right \rangle =-\left \langle Ay,x\right \rangle $ for all
$x,y\in V$, or equivalently, $\left \langle Ax,x\right \rangle =0$ for all $x\in
V$.

Given a set $S\subseteq X$, $\overline{S}$ will denote its closure.

\section{The main result}

\begin{theorem}
	\label{Th: basic}Let $X$ be a normed space and $T:X\rightrightarrows X^{\ast}$
	be an operator. Given $\left(  u,u^{\ast}\right)  \in \operatorname*{gph}T$,
	let $T^{\prime}$ be the translation $T^{\prime}\left(  x\right)  :=T\left(
	x+u\right)  -u^{\ast}$, $x\in X$ and $\hat{T}$ be the reduction of $T^{\prime
	}$ on the space $V=\left[  \operatorname*{dom}T^{\prime}\right]  $. The
	operator $T$ is bimonotone iff there exists a single-valued, skew symmetric
	operator $A:V\rightarrow V^{\ast}$ such that $\hat{T}\left(  x\right)  =Ax$
	for all $x\in \operatorname*{dom}T^{\prime}$.
\end{theorem}

\begin{proof}
	Assume that $T$ is bimonotone. Obviously, the operator $T^{\prime}$ is
	bimonotone and $0\in T^{\prime}\left(  0\right)  $. Also, $\operatorname*{dom}
	T^{\prime}=\operatorname*{dom}T-u$.
	
	Set $V=\left[  \operatorname*{dom}T^{\prime}\right]  $. Let $\hat
	{T}:V\rightrightarrows V^{\ast}$ be the reduction of $T^{\prime}$.
	
	Let us check that $\hat{T}$ is bimonotone, i.e.,
	\begin{equation}
		\forall \left(  x,x^{\ast}\right)  ,\left(  y,y^{\ast}\right)  \in
		\operatorname*{gph}\hat{T},\quad \left \langle x^{\ast}-y^{\ast}
		,x-y\right \rangle =0. \label{rel:basic}
	\end{equation}

	Indeed, if $\left(  x,x^{\ast}\right)  $, $\left(  y,y^{\ast}\right)
	\in \operatorname*{gph}\hat{T}$, then $x^{\ast},y^{\ast}$ can be written as
	$x^{\ast}=J^{\ast}x_{1}^{\ast}$, $y^{\ast}=J^{\ast}y_{1}^{\ast}$ with
	$x_{1}^{\ast}\in T^{\prime}\left(  Jx\right)  $, $y_{1}^{\ast}\in T^{\prime
	}\left(  Jy\right)  $. Hence,
	\[
	\left \langle x^{\ast}-y^{\ast},x-y\right \rangle =\left \langle J^{\ast}
	x_{1}^{\ast}-J^{\ast}y_{1}^{\ast},x-y\right \rangle =\left \langle x_{1}^{\ast
	}-y_{1}^{\ast},Jx-Jy\right \rangle =0
	\]
	so $\hat{T}$ is also bimonotone.
	
	Since $\left(  0,0\right)  \in \operatorname*{gph}\hat{T}$, it follows from
	(\ref{rel:basic}) that $\left \langle x^{\ast},x\right \rangle =0$ for all
	$\left(  x,x^{\ast}\right)  \in \operatorname*{gph}\hat{T}$. Then
	(\ref{rel:basic}) implies that for all $\left(  x,x^{\ast}\right)  $ and
	$\left(  y,y^{\ast}\right)  $ in $\operatorname*{gph}\hat{T}$,
	\begin{equation}
		\left \langle x^{\ast},y\right \rangle =-\left \langle y^{\ast},x\right \rangle .
		\label{rel:basic'}
	\end{equation}

	Let us show that for all $x\in \operatorname*{dom}\hat{T}$, $\hat{T}\left(
	x\right)  $ is a singleton. Indeed, pick any $\left(  y,y^{\ast}\right)
	\in \operatorname*{gph}\hat{T}$. If $x_{1}^{\ast},x_{2}^{\ast}\in \hat{T}\left(
	x\right)  $, then replacing successively $x^{\ast}$ by $x_{1}^{\ast}$,
	$x_{2}^{\ast}$ in (\ref{rel:basic'}) entails
	\[
	\left \langle x_{1}^{\ast}-x_{2}^{\ast},y\right \rangle =0.
	\]

	This is true for all $y\in \operatorname*{dom}\hat{T}$, so it is true for all
	$y\in \left[  \operatorname*{dom}T^{\prime}\right]  =V$. Thus, $\hat{T}\left(
	x\right)  $ is a singleton for all $x\in \operatorname*{dom}\hat{T}$ and
	$\hat{T}$ can be considered as a single-valued operator. From now on, $\hat
	{T}\left(  x\right)  $ will denote the element in the image of $x$.
	
	Now we define a linear operator $A:V\rightarrow V^{\ast}$ as follows: For
	every $x\in V$ we write $x=\sum_{i\in I}\alpha_{i}x_{i}$ ($I$ finite) with
	$\alpha_{i}\in \mathbb{R}$, $x_{i}\in \operatorname*{dom}\hat{T}$, and set
	$Ax=\sum_{i\in I}\alpha_{i}\hat{T}\left(  x_{i}\right)  $. The operator is
	well-defined, because if $x=\sum_{j\in J}\beta_{j}y_{j}$ with $\beta_{j}
	\in \mathbb{R}$, $y_{j}\in \operatorname*{dom}\hat{T}$ ($J$ finite) is another
	expansion of $x$, then $0=\sum_{i\in I}\alpha_{i}x_{i}-\sum_{j\in J}\beta
	_{j}y_{j}$. For every $z\in \operatorname*{dom}\hat{T}$, (\ref{rel:basic})
	implies
	\begin{align*}
		\left \langle \sum_{i\in I}\alpha_{i}\hat{T}\left(  x_{i}\right)  -\sum_{j\in
			J}\beta_{j}\hat{T}\left(  y_{j}\right)  ,z\right \rangle  &  =\sum_{i\in
			I}\alpha_{i}\left \langle \hat{T}\left(  x_{i}\right)  ,z\right \rangle
		-\sum_{j\in J}\beta_{j}\left \langle \hat{T}\left(  y_{j}\right)
		,z\right \rangle \\
		&  =-\sum_{i\in I}\alpha_{i}\left \langle \hat{T}\left(  z\right)
		,x_{i}\right \rangle +\sum_{j\in J}\beta_{j}\left \langle \hat{T}\left(
		z\right)  ,y_{j}\right \rangle \\
		&  =0.
	\end{align*}

	This holds for every $z\in \operatorname*{dom}\hat{T}$ and eventually for every
	$z\in V$. Hence, $\sum_{i\in I}\alpha_{i}\hat{T}\left(  x_{i}\right)
	=\sum_{j\in J}\beta_{j}\hat{T}\left(  y_{j}\right)  $, i.e., the value of $Ax$
	is independent of the expansion of $x$. Note that for every $x\in
	\operatorname*{dom}\hat{T}$, $x$ has the expansion $x=x$ so $Ax=\hat{T}\left(
	x\right)  $.
	
	It can be readily checked that $A$ is linear. Let us check that it is skew
	symmetric: For $x,y\in V$ with the expansions $x=\sum_{i\in I}\alpha_{i}x_{i}
	$, $y=\sum_{j\in J}\beta_{j}y_{j}$ ($x_{i},y_{j}$ in $\operatorname*{dom}
	\hat{T}$), using (\ref{rel:basic}) we find:
	\begin{align*}
		\left \langle Ax,y\right \rangle  &  =\sum_{i\in I}\sum_{j\in J}\alpha_{i}
		\beta_{j}\left \langle Ax_{i},y_{j}\right \rangle =\sum_{i\in I}\sum_{j\in
			J}\alpha_{i}\beta_{j}\left \langle \hat{T}\left(  x_{i}\right)  ,y_{j}
		\right \rangle \\
		&  =-\sum_{i\in I}\sum_{j\in J}\alpha_{i}\beta_{j}\left \langle \hat{T}\left(
		y_{j}\right)  ,x_{i}\right \rangle =-\left \langle Ay,x\right \rangle .
	\end{align*}

	Thus, one direction of the proposition is proved. The converse is easy: Let
	$T$ be an operator such that the operator $\hat{T}$ defined as above satisfies
	$\hat{T}\left(  x\right)  =Ax$ on $V$, where $A$ is single-valued, linear and
	skew symmetric. Two arbitrary elements of $\operatorname*{gph}T^{\prime}$ can
	be written as $\left(  Jx,x^{\ast}\right)  $ and $\left(  Jy,y^{\ast}\right)
	$ with $x,y\in V$. We have:
	\begin{align*}
		\left \langle x^{\ast}-y^{\ast},Jx-Jy\right \rangle  &  =\left \langle J^{\ast
		}x^{\ast}-J^{\ast}y^{\ast},x-y\right \rangle =\left \langle \hat{T}\left(
		x\right)  -\hat{T}\left(  y\right)  ,x-y\right \rangle \\
		&  =\left \langle Ax-Ay,x-y\right \rangle =0.
	\end{align*}

	Hence $T^{\prime}$ is bimonotone, so $T$ is also bimonotone.
\end{proof}

Remarks:

\begin{enumerate}
	\item The elements of $T^{\prime}\left(  x\right)  $ are extensions to the
	whole space $X$ of the unique element $Ax$ (seen as a continuous linear
	functional on $V$).
	
	\item Assume that $X$ is a Banach space. Then also $\overline{V}$ is a Banach
	space. The operator $A$ can be extended to a skew symmetric operator on
	$\overline{V}$, if and only if it is continuous. Indeed, a skew symmetric
	operator defined on the Banach space $\overline{V}$ is continuous, as a result
	of the closed graph theorem. To see this, we just copy the well-known proof
	for symmetric operators: If $\left(  x_{n}\right)  $ is a sequence in
	$\overline{V}$ such that $x_{n}\rightarrow x$ and $Ax_{n}\rightarrow x^{\ast}
	$, then for every $y\in \overline{V}$, $\left \langle x^{\ast},y\right \rangle
	=\lim \left \langle Ax_{n},y\right \rangle =-\lim \left \langle Ay,x_{n}
	\right \rangle =-\left \langle Ay,x\right \rangle =\left \langle Ax,y\right \rangle
	$. Thus, $x^{\ast}=Ax$. Hence $A$ has a closed graph and by the closed graph
	theorem, $A$ is continuous. A skew symmetric operator defined on a normed
	space which is not complete can be discontinuous, as shown by the following
	example: Take $X=L^{2}\left(  0,1\right)  $ with the $\left \Vert
	\cdot \right \Vert _{2}$ norm, $V=\mathcal{C}_{c}^{1}\left(  0,1\right)  $ (the
	space of compactly supported $\mathcal{C}^{1}$ functions on $\left(
	0,1\right)  $), $A=\frac{d}{dx}$ defined on $V$. Here, $\overline{V}=X$, but
	$A$ is not continuous so it cannot be extended to $X$.
\end{enumerate}

When the linear span of the domain of a bimonotone operator $T$ is dense in
$X$, then $T$ has a very simple form:

\begin{corollary}
	\label{cor:basic}Let $X$ be a normed space and $T:X\rightrightarrows X^{\ast}$
	be an operator such that $V:=\left[  \operatorname*{dom}T\right]  $ is dense
	in $X$. The operator $T$ is bimonotone iff there exists a single-valued, skew
	symmetric operator $A:X\rightarrow X^{\ast}$ with domain $V$ and $v^{\ast}\in
	X^{\ast}$ such that $T\left(  x\right)  =Ax+v^{\ast}$ for all $x\in
	\operatorname*{dom}T$.
\end{corollary}

\begin{proof}
	Since the subspace $V$ is dense in $X$, we can identify $V^{\ast}$ with
	$X^{\ast}$. Take any $\left(  u,u^{\ast}\right)  \in \operatorname*{gph}T$ and
	define as before $T^{\prime}:X\rightrightarrows X^{\ast}$ by $T^{\prime
	}\left(  x\right)  =T\left(  x+u\right)  -u^{\ast}$, $x\in X$. Since $\left[
	\operatorname*{dom}T^{\prime}\right]  =V$ is dense in $X$, it is clear that
	$\tilde{T}=T^{\prime}$. By Theorem \ref{Th: basic}, there exists a
	single-valued, skew symmetric operator $A:X\rightarrow X^{\ast}$ with domain
	$V$ such that $T^{\prime}\left(  x\right)  =Ax$ for all $x\in
	\operatorname*{dom}T^{\prime}$. This means that in particular, $T$ and
	$T^{\prime}$ are single-valued. For every $x\in \operatorname*{dom}T$,
	$T\left(  x\right)  =T^{\prime}\left(  x-u\right)  +u^{\ast}=Ax-Au+u^{\ast}$.
	Setting $v^{\ast}=u^{\ast}-Au$ we obtain the result.
\end{proof}

As an immediate application, we obtain a generalization of \cite[Prop. 3.6
(i)(ii)]{Martinez}:

\begin{proposition}
	Let $T:X\rightrightarrows X^{\ast}$ be paramonotone and bimonotone. Then, we have
	
	(i) If $\left[  \operatorname*{dom}T\right]  $ is dense in $X$, then
	$\operatorname*{range}T$ is a singleton;
	
	(ii) If $\left[  \operatorname*{range}T\right]  $ is dense in $X^{\ast}$, then
	$\operatorname*{dom}T$ is a singleton.
\end{proposition}

\begin{proof}
	(i) If $T$ is paramonotone and bimonotone, then it is immediate from the
	definitions that it is constant on its domain (see also \cite[Prop.
	3.5]{Martinez}). According to Corollary \ref{cor:basic} the operator $T$ is
	single-valued, so the assertion follows.
	
	(ii) It is obvious that $T$ is paramonotone and bimonotone iff its inverse
	$T^{-1}:X^{\ast}\rightrightarrows X$ is paramonotone and bimonotone (where we
	consider $X$ as a subset of $X^{\ast \ast}$). Since $\operatorname*{range}
	T=\operatorname*{dom}T^{-1}$ and $\operatorname*{dom}T=\operatorname*{range}
	T^{-1}$, (ii) follows from (i).
\end{proof}

\end{document}